\documentstyle[12pt]{amsart}

\def\sw#1{{\sb{(#1)}}}

\def\<{{\langle}}
\def\>{{\rangle}}

\def\eps{\epsilon}

\def\note#1{{}}

\def\note#1{}
\def\M{{\cal M}}

\def\cC{{\cal C}}

\def\End{{\rm End}}

\def\Label{\label}
\begin{document}

\baselineskip 22pt

\newtheorem{proposition}{Proposition}[section]
\newtheorem{lemma}[proposition]{Lemma}
\newtheorem{corollary}[proposition]{Corollary}
\newtheorem{theorem}[proposition]{Theorem}

\theoremstyle{definition}
\newtheorem{definition}[proposition]{Definition}
\newtheorem{example}[proposition]{Example}

\theoremstyle{remark}
\newtheorem{remark}[proposition]{Remark}

\newcommand{\Section}{\setcounter{definition}{0}\section}
\newcounter{c}
\renewcommand{\[}{\setcounter{c}{1}$$}
\newcommand{\etyk}[1]{\vspace{-7.4mm}$$\begin{equation}\Label{#1}
\addtocounter{c}{1}}
\renewcommand{\]}{\ifnum \value{c}=1 $$\else \end{equation}\fi}

\title[Bialgebroids]{Bialgebroids, $\times_{A}$-bialgebras
and duality}
\author{Tomasz Brzezi\'nski}
\address{Department of Mathematics, University of Wales Swansea,
Singleton Park, Swansea SA2 8PP, U.K.}
\email{T.Brzezinski@@swansea.ac.uk}
\author{Gigel Militaru}
\address{Faculty of Mathematics, University of Bucharest, Str.
Academiei 14, RO-70109 Bucharest 1, Romania}
\email{gmilit@@al.math.unibuc.ro}
\subjclass{16W30, 13B02}
\begin{abstract}
An equivalence between Lu's bialgebroids, Xu's bialgebroids with
an anchor and Takeuchi's $\times_{A}$-bialgebras is explicitly
proven. A new class of examples of bialgebroids is constructed. A
(formal) dual of a bialgebroid, termed bicoalgebroid, is defined.  A
weak Hopf algebra is shown to be an example of such a bicoalgebroid.
\end{abstract}
\maketitle
\section{Introduction}
For some time various generalisations of the notion of a bialgebra,
in which a bialgebra is required to be a bimodule but not
necessarily an algebra over a (noncommutative) ring
have been considered. Motivated by the problem of classification of
algebras, a definition of a generalised Hopf algebra  was first
proposed by Sweedler \cite{Swe:alg} and later generalised by
Takeuchi \cite{Tak:gro}. This was based on a new definition of a tensor
product over noncommutative rings, termed the $\times_{A}$-product.
Several years later, motivated by some problems in algebraic topology
Ravenel introduced the notion of a commutative Hopf algebroid
\cite{Rav:com}, which is a special case of the Takeuchi construction.
With the growing interest in quantum groups, bialgebroids were
discussed  in the context of noncommutative \cite{Mal},
and Poisson geometry. In the latter case,
the most general definitions were given by Lu \cite{Lu} and Xu
\cite{Xu}. Another generalisations of finite Hopf algebras, termed
weak Hopf algebras, appeared in relation to integrable spin chains and
classification of subfactors of von Neumann algebras \cite{BohSzl:coa}
\cite{Nil:axi}.
In \cite{EtiNik:dyn} weak Hopf algebras have been
shown to be examples of Lu's bialgebroids.

The aim of the present paper is threefold. Since there is a number of
different definitions of generalised bialgebras, it is important to
study what are the relations between these definitions. Thus our
first aim (Section~2) is to collect these different
definitions and make it clear that
the notions of a Takeuchi's
$\times_{A}$-bialgebra,
Lu's bialgebroid, and  Xu's bialgebroid with
an anchor are equivalent to each other (Section~3).
Although this fact in itself
seems to be not new (cf.\ \cite[p.\ 273]{Sch:dua}, where the
equivalence of the first and second notions is attributed to P.\ Xu
\cite{Xu}), to the best of our
knowledge, there is no explicit and complete proof of this
equivalence in the literature.  Hereby we provide such a proof,
and hope that this clarifies some minor misunderstandings
in the field (e.g.\ it seems to be
claimed in \cite[p.\ 546]{Xu} that Lu's bialgebroid is equivalent
to Xu's bialgebroid {\em without an anchor}).
Our second aim (Section~4) is to
construct new examples of bialgebroids. We show how to associate a
bialgebroid to a braided commutative algebra in the category of
Yetter-Drinfeld modules. This result generalises an example
considered by Lu. In fact we show that the smash product of a Hopf
algebra with an algebra in the Yetter-Drinfeld category is a
bialgebroid if and only if the algebra is braided commutative.
In particular, a bialgebroid
over braided symmetric algebra $S_R (n)$ is associated
to any solution of the quantum Yang-Baxter equation
$R\in M_n(k) \otimes M_n(k)$.
Our third aim (Section~5) is to propose a notion that is
dual to a bialgebroid. We
term such an object a {\em bicoalgebroid}. It is
well-known that a bialgebra is a self-dual notion in the following
sense. The axioms
of a bialgebra are invariant under formal reversing of the arrows in
the commutative diagrams that constitute the definition of a
bialgebra. In the case of a bialgebroid such a formal operation on
commutative diagrams produces a new object. We belive that this object
will play an important role in constructing a self-dual generalisation
of a bialgebra which should involve both a bialgebroid and a
bicoalgebroid.

\section{Preliminaries}
\subsection{Notations}\label{subs.notation}
All rings in this
paper have 1, a ring map is assumed to respect 1, and all
modules over a ring are assumed to be unital.
For a ring $A$,  $\M_A$ (resp.\ ${}_A\M$) denotes
the category of right (resp.\ left) $A$-modules, and ${}_A\M_{A}$ 
denotes the category of $(A,A)$-bimodules.
The action of $A$ is denoted by a dot between elements.

Throughout the paper $k$ denotes a commutative ring.
Unadorned tensor product
is over $k$. For a $k$-algebra $A$ we use $m_A$ to denote the
product as a map and $1_A$ to denote unit both as an element of
$A$ and as a map $k\to A$, $\alpha\to \alpha 1$. ${\rm End}(A)$
denotes the algebra of $k$-linear endomorphisms of $A$.
For a
$k$-coalgebra $C$ we use $\Delta$ to denote the coproduct,
$\eps$ to denote the counit;
$\M^C$ will be the category of right $C$-comodules.
We use the Sweedler notation,
i.e. $\Delta(c) = c\sw 1\otimes c\sw 2$ for coproducts, and 
$\rho^M(m) = m_{<0>} \otimes m_{<1>}$ for coactions (summation 
understood).

Let $A$ be a $k$-algebra. Recall from \cite{Swe:pre} that an
{\em $A$-coring} is an $(A,A)$-bimodule
$\cC$ together with $(A,A)$-bimodule maps $\Delta_\cC:\cC\to
\cC\otimes_A\cC$ called a coproduct and $\eps_\cC:\cC\to A$ called a
counit, such that
$$
(\Delta_\cC\otimes_A\cC)\circ\Delta_\cC = (\cC\otimes_A\Delta_\cC)\circ
\Delta_\cC, \quad (\eps_\cC\otimes_A\cC)\circ\Delta_\cC =
(\cC\otimes_A\eps_\cC)\circ \Delta_\cC = \cC.
$$
Let $R$ be a $k$-algebra. Recall from \cite{Swe:alg}, \cite{Tak:gro}
that an {\em $R$-ring} or an {\em algebra over $R$} is a pair
$(U, i)$, where $U$ is a $k$-algebra  and $i: R\to U$ is an algebra map.
If $(U, i)$ is an $R$-ring then $U$ is an $(R,R)$-bimodule  with the
structure provided by the map $i$,
$r\cdot u \cdot r':= i(r)ui(r')$. A map  of $R$-rings $f: (U, i)\to (V,j)$
is a
$k$-algebra map $f:U\to V$ such that $f\circ i=j$. Equivalently, a map
of $R$-rings is a $k$-algebra map that is a left or right $R$-module
map. Indeed, clearly if $f:(U,i)\to (V,j)$ is a map of $R$-rings it is
an algebra and
$R$-bimodule map. Conversely, if $f$ is a left $R$-linear algebra
map then for all $r\in R$, $f(i(r)) = f(r\cdot 1_U) = r\cdot f(1_U) =
j(r)$, and similarly in the right $R$-linear case.

\subsection{Algebras over enveloping algebras: $A^{e}$-rings.}\label{subs.ae}
Let $A$ be an algebra and $\bar{A}=A^{op}$ the opposite algebra.
For $a\in A$, $\bar{a}\in \bar{A}$ is the same $a$ but now viewed as an
element in $\bar{A}$, i.e. $a\to \bar{a}$ is an antiisomorphism of
algebras. Let
$
A^{e}=A\otimes \bar{A}
$
be the enveloping algebra of $A$.
Note that a pair $(H,i)$ is  an $A^{e}$-ring if and
only if there exist an algebra map $s:A\to H$ and an anti-algebra map
$t: A\to H$,
such that
$s(a)t(b)=t(b) s(a)$,
for all $a,b\in A$. Explicitly, $s(a) = i(a\otimes 1)$ and $t(a) =
i(1\otimes \bar{a})$, and, conversely,  $i(a\otimes \bar{b})=s(a)t(b)$.
This simple observation
shows that the bialgebroids  of Lu \cite{Lu} (cf.\
Section~\ref{subs.bialgebroids} below)
and $\times_A$-bialgebras of Takeuchi \cite{Tak:gro} (cf.\
Section~\ref{subs.x-bialg} below)
have the same input data.

In the sequel, the expression ``let $(H, s, t)$ be an $A^{e}$-ring" will
be understood  to mean an algebra $H$ with algebra maps $s,t:A \to H$
 as described above.
$A$ is called a
{\em base algebra}, $H$ a {\em total algebra}, $s$ the {\em source}
map and $t$ the {\em target} map.

A standard example of an $A^e$-ring is provided by $\End(A)$. In this case
$i:A\otimes \bar{A}\to \End(A)$, $i(a\otimes \bar{b})(x)=axb$.
The source and the target come out as, $s(a)(x)=ax$ and
$t(b)(x)=xb$.
It follows that $\End(A)$ is an $A^{e}$-bimodule via $i$.
In particular,
$\End(A)$ is a left $A^e$-module via
\begin{equation}\label{8}
(a\cdot f)(b)=af(b), \qquad
(\bar{a}\cdot f)(b)=f(b)\bar{a},
\end{equation}
for all $a, b\in A$ and $f\in \End(A)$.

As explained at the end of Section~\ref{subs.notation},
an $A^{e}$-ring $(H, s, t)$  is an
$A^{e}$-bimodule. This means that $H$ is an $A$-bimodule with explicit
actions, $a\cdot h \cdot b = s(a)h s(b)$, and also that it is an
$\bar{A}$-bimodule with actions
$\bar{a}\cdot h \cdot \bar{b}= t(a)h t(b)$,
for all $a,b\in A$ and $h\in H$.

\subsection{The key $A^{e}$-ring associated to an $A^{e}$-ring}
Let $(H, s, t)$  be an $A^{e}$-ring and view $H$ as an
$A$-bimodule, with the left $A$-action given by the source map $s$,
and the right $A$-action
which descends from the left $\bar{A}$-action given by the target map
$t$, i.e.,
\begin{equation}\label{16'}
a\cdot h = s(a) h, \qquad h\cdot a = t(a)h, \qquad \forall a\in A, h\in H.
\end{equation}
Consider an Abelian group $H\otimes_A H$ which is an
$A$-bimodule via the following actions
$a\cdot ' (g\otimes_A h) = g t(a) \otimes_A h$, and
$(g\otimes_A h) \cdot ' a = g \otimes_A hs(a)$,
for all $a\in A$ and $g\otimes_A h\in H\otimes_A H$. Define
$\Gamma = \Gamma (H, s, t) := (H\otimes_A H)^A$, i.e.,
$$
\Gamma = \{ g =
\sum g^{1}\otimes_{A} g^{2} \in H\otimes_A H\; |\; \forall a\in A ,\;
\sum g^1 t(a) \otimes_A g^2 = \sum g^1\otimes_A g^2 s(a)\}.
$$
The next proposition can be proved directly (cf.\
  \cite[Proposition 3.1]{Tak:gro}).

\begin{proposition}\label{2.1}
Let $(H, s,t)$ be an $A^{e}$-ring. Then  $\Gamma = \Gamma (H, s, t)$ is an
$A^{e}$-ring with the algebra structure
$
(\sum g^1\otimes_A g^2) (\sum h^1\otimes_A h^2)=
\sum g^1h^1\otimes_A g^2h^2,
$
the unit $1_H\otimes_A 1_H$ and the algebra map
$i: A\otimes \bar{A} \to \Gamma$,
$a \otimes_A \bar{b}\mapsto s(a)\otimes_A t(b)$.
\end{proposition}

\subsection{Bialgebroids}\label{subs.bialgebroids}
Let $(H, s, t)$ be an $A^{e}$-ring and view $H \in {}_A\M_A$ using the
actions in equations (\ref{16'}).  Also, view $H\otimes_A H \in {}_A\M_A$
with the natural actions $a \cdot (g\otimes_A h) \cdot b= s(a) g\otimes_A
t(b)h$.

\begin{definition}[\cite{Lu}]\label{def.Lu}
Let $(H, s, t)$ be an $A^{e}$-ring. We say that $(H, s, t, \Delta, \eps)$ is
an {\em $A$-bialgebroid} iff

(B1)~~~~~ $(H,\Delta, \eps)$ is an $A$-coring;

(B2)~~~~~ ${\rm Im}(\Delta) \subseteq \Gamma (H, s, t)$ and the corestriction
of the coproduct $\Delta :H \to \Gamma (H, s, t)$ is an algebra map;

(B3)~~~~~$\eps (1_H)=1_A$ and for all $g$, $h\in H$
\begin{equation}\label{25}
\eps(gh) =\eps \Bigl( g s( \eps (h) ) \Bigl) =
\eps \Bigl( g t( \eps (h) ) \Bigl).
\end{equation}

An {\em antipode} for an $A$-bialgebroid $H$ is an anti-algebra map
$\tau :H \to H$ such that

(ANT1)~~~~~ $\tau \circ t =s$;

(ANT2)~~~~~$m_H \circ (\tau \otimes H) \circ \Delta = t\circ \eps
\circ \tau$;

(ANT3)~~~~~There exists a section $\gamma :H\otimes_A H \to H\otimes H$ of
the natural
projection  $H\otimes H \to H\otimes_A H$ such that
$m_H \circ (H\otimes \tau)\circ \gamma \circ \Delta = s\circ \eps$.
An $A$-bialgebroid with an antipode is called a {\em Hopf algebroid}.
\end{definition}
The notion of a bialgebroid was introduced by  J.-H.\ Lu in \cite{Lu}.
Condition
(B2), in its present form, was first stated by P.\ Xu \cite{Xu},
while (B3) in this form appears in \cite{Szl} (in a slightly
different convention though).
\begin{remark}
1) As explained in Section~\ref{subs.notation}
the axiom (B1) requires that $\Delta :H \to H\otimes_A H$ and
$\eps :H \to A$ are maps in ${}_A\M_A$,
$\Delta$ is coassociative
and $\eps$ is a counit for $\Delta$. The counit property
explicitly means that for all $h\in H$,
\begin{equation}\label{24}
s \Bigl(\eps (h_{(1)}) \Bigl) h_{(2)} = t\Bigl(\eps (h_{(2)})
\Bigl) h_{(1)} = h.
\end{equation}
The first condition of (B2) explicitly means that
$h_{(1)} t(a)\otimes_A h_{(2)} = h_{(1)} \otimes_A h_{(2)} s(a),$
for all $a\in A$ and $h\in H$.

2) Definition~\ref{def.Lu} is equivalent to   \cite[Definition~2.1]{Lu}.
Indeed, in \cite[Proposition 3.2]{Xu} it is shown that (B2)
in Definition~\ref{def.Lu} is
equivalent to condition 4 in \cite[Definition~2.1]{Lu}, which states
that the kernel of
the map
$\Phi : H\otimes H\otimes H \to H\otimes_A H$,
$\Phi(g \otimes h \otimes l)=\Delta (g) (h\otimes l)$
is a left ideal of $H\otimes \bar{H} \otimes \bar{H}$.

On the other hand condition (\ref{25}) of (B3) is equivalent to
the condition in \cite{Lu} that $\ker(\eps)$ is a left ideal in $H$.
Indeed, suppose that $\ker(\eps)$ is a left ideal in $H$. Then
using that $\eps$ is an $A$-bimodule map and
that $\eps(1)=1$ one obtains that
$h- s(\eps (h)), h- t(\eps (h)) \in \ker(\eps)$,
so (\ref{25}) holds. The converse is obvious.

3) The facts that $\eps$ preserves unit and is an $A$-bimodule map
imply that $s$ and $t$ are sections of $\eps$, i.e.,
$\eps(s(a)) = \eps(t(a)) = a$ for all $a\in A$. Using
this fact and (\ref{25}) one easily finds that
$\eps(g s(a)) = \eps(g t(a))$,
for all $a\in A$ and $g\in H$
Similarly, the facts that $\Delta$ is a unital and $A$-bimodule map
imply that
\begin{equation}\label{4.2}
    \Delta (s(a)) = s(a) \otimes_A 1_H, \qquad
    \Delta (t(a)) = 1_H \otimes_A t(a).
\end{equation}

4) For an $A^{e}$-ring $(H, s, t)$, let $F: {}_H\M \to {}_A\M_A$,
be the
restriction of scalars functor. The actions of $A$ on $H$ are  given
by  equations  (\ref{16'}).
If $(H, s, t, \Delta, \eps)$ is an $A$-bialgebroid, then
${}_H\M$ has a monoidal structure such that $F$ is a strict monoidal
functor. For all $M, N\in {}_H\M$, the tensor product $M\otimes_A N$ is
in ${}_H\M$ via
$h\cdot (m\otimes_A n) = h_{(1)}\cdot m\otimes_A h_{(2)}\cdot n$. The
 right hand side is well defined because ${\rm Im}(\Delta) \subseteq
 \Gamma$.  $A$ is the unit object, when viewed in
 ${}_H\M$ via the action
$h\triangleright a = \eps (h s(a))= \eps (h t(a)),$
for all $h\in H$, $a\in A$. The fact that this is an action follows from
(\ref{25}) (cf.\ \cite[Eq.~(8)]{Lu}).
\end{remark}

Next we recall the notion  of a bialgebroid with an anchor from \cite{Xu}.
For a $k$-algebra $A$, $\End(A)$  is an $A^{e}$-ring
as described in Section~\ref{subs.ae}. In fact,
we are interested only in the
structure
of $\End(A)$ as a left ${A\otimes \bar{A}}$-module given by equations
(\ref{8}).  When, following  \cite{Xu},  $\End(A)$ is viewed as an
$A$-bimodule, the structure maps are
\begin{equation}\label{9'}
(a\cdot f)(b)=af(b), \quad
(f\cdot a )(b)=f(b)a,
\end{equation}
for all $a$, $b\in A$ and $f\in \End(A)$. As before the total
algebra $H$ of an $A^{e}$-ring $(H, s, t)$, is viewed as an $A$-bimodule
by  equations (\ref{16'}).

Leaving aside the redundant part of \cite[Definition 3.4]{Xu} that
arises from condition ii) there, \cite[Proposition 3.3]{Xu}, and
from the fact that a counit of a coring is unique the following
definition is the same as \cite[Definition 3.4]{Xu}

\begin{definition}\label{2.4}
Let $(H, s, t)$ be an $A^{e}$-ring. $(H, s, t, \Delta, \mu)$ is called
an {\em $A$-bialgebroid with an anchor $\mu$} if

(BA1)~~~~~ $\Delta :H \to H\otimes_A H$ is a coassociative  $A$-bimodule map;

(BA2)~~~~~ ${\rm Im}(\Delta) \subseteq \Gamma (H, s, t)$ and its corestriction
$\Delta :H \to \Gamma (H, s, t)$ is an algebra map;

(BA3)~~~~~$\mu :H\to \End(A)$ is an algebra and an $A$-bimodule map such that

~~~~~~~~~~~~~~~~~~~~~~~(A1)~~~~~$s(h_{(1)}\triangleright a)h_{(2)}=hs(a)$,

~~~~~~~~~~~~~~~~~~~~~~~(A2)~~~~~$t(h_{(2)}\triangleright a)h_{(1)}=ht(a)$,

where $\mu(h)(a) = h\triangleright a$, for all $a\in A$ and $h\in H$.
\end{definition}

Note that the left hand sides of (A1) and (A2) are well defined since
$\mu$ is an $A$-bimodule map, i.e., for all $a$, $b\in A$ and $h\in H$
\begin{equation}\label{29}
\mu (s(a) h)(b) = a\mu(h)(b), \qquad
\mu (t(a) h)(b) = \mu(h)(b)a.
\end{equation}

\subsection{$\times_A$-bialgebras}\label{subs.x-bialg}
The notion of a $\times_A$-bialgebra was first introduced
by M.E.\ Sweedler \cite{Swe:alg} for a commutative $A$ and then
generalised by M.\ Takeuchi
\cite{Tak:gro} to an arbitrary $A$. In this section we briefly recall
Takeuchi's definition (see \cite{Tak:gro} for details).

Let $M$ and $N$ be  $A^{e}$-bimodules. Following MacLane, let
$$
\int_a {}_{\bar{a}}M \otimes {}_aN :=
M\otimes N /<\{\bar{a} m\otimes n -m \otimes an | \forall a\in A \}>,
$$
$$
\int^b M_{\bar{b}} \otimes N_b :=
\{\sum_{i} m_i\otimes n_i \in M\otimes N | \forall b\in A, \;
\sum_{i} m_i\bar{b} \otimes n_i = \sum m_i \otimes n_i b \}.
$$
Define,
$$
M\times_A N := \int^b \int_a {}_{\bar{a}}M_{\bar{b}} \otimes {}_aN_b.
$$
The operation
$-\times_A - : {}_{A^{e}}\M_{A^{e}} \times {}_{A^{e}}\M_{A^{e}} \to
{}_{A^{e}}\M_{A^{e}}$ is a bifunctor.
Here, for $M, N\in {}_{A^{e}}\M_{A^{e}}$, the product
$M\times_A N$ is  in ${}_{A^{e}}\M_{A^{e}}$ with the actions given by
\begin{equation}\label{32}
(a\otimes\bar{a})\cdot (\sum_{i} m_i \otimes_A n_i)\cdot (b\otimes
\bar{b}) = \sum_{i} am_i b \otimes_A \bar{a}n_i\bar{b}.
\end{equation}
For any two $A^e$-rings $(U, i)$, $(V, j)$, $U\times_A V$ is an $A^e$-ring via
the well defined algebra map
$A\otimes \bar{A} \to U\times_A V$, $a\otimes \bar{b}\to i(a) \otimes_A
j(\bar{b})$. Note that if $(H, s, t)$ is an $A^e$-ring and $H$ is considered as
an $A^e$-bimodule via the
actions described at the end of Section~\ref{subs.ae}, then
$H\times_A H = \Gamma (H, s, t)$.

For $M$, $N$ and $P \in {}_{A^{e}}\M_{A^{e}}$ define
$$
M\times_A P\times_A N := \int^{s,u}\int_{r,t}
{}_{\bar{r}}M_{\bar{s}} \otimes {}_{r, \bar{t}}P_{s, \bar{u}}
\otimes {}_tN_u.
$$
There exist natural maps
$$
\alpha : (M\times_A P)\times_A N \to M\times_A P\times_A N , \quad
\alpha' :M\times_A (P\times_A N) \to M\times_A P\times_A N.
$$
The maps $\alpha$, $\alpha'$ are not isomorphisms in general.
Since $\End(A)$ is an $A^e$-ring (cf.\ Section~\ref{subs.ae}),
it is an $A^e$-bimodule, so
one can define the maps
$$
\theta : M\times_A \End(A) \to M, \quad
\theta (\sum_{i} m_i \otimes_A f_i) = \sum_{i} \overline{f_i(1)} m_i,
$$
$$
\theta' : \End(A) \times_A M \to M, \quad
\theta' (\sum_{i} f_i \otimes_A m_i)= \sum_{i} f_i(1) m_i,
$$
Following \cite{Tak:gro} a triple $(L, \Delta, \mu)$ is called a {\em
$\times_A$-coalgebra} iff
$L$ is an $A^e$-bimodule and
$\Delta :L \to L\times_A L$, $\mu: L\to \End(A)$
are $A^e$-bimodule maps such that
$$
\alpha \circ (\Delta \times_A L) \circ \Delta =\alpha' \circ (L\times_A
\Delta) \circ \Delta,
\quad
\theta \circ (L\times_A \mu) \circ \Delta = L = \theta' \circ (\mu\times_A
L) \circ \Delta.
$$

\begin{remark}[\cite{Sch:dua}]\label{2.7}
Let $L$ be an $A^e$-bimodule and
$\Delta :L \to L\times_A L$ and $\mu: L\to \End(A)$
$A^e$-bimodules maps. Let $i: L\times_A L\to L\otimes_A L$ be the
canonical inclusion. Then $(L, \Delta, \mu)$ is a
$\times_A$-coalgebra if and only if
$(L, \Delta', \eps_{\mu})$ is an $A$-coring, where
$\Delta' = i\circ\Delta$ and
$\eps_{\mu} (l) = \mu(l)(1_A)$ .
\end{remark}

\begin{definition}[\cite{Tak:gro}]\label{2.8}
Let $(H, s, t)$ be an $A^{e}$-ring. $(H, s, t, \Delta, \mu)$ is called
a {\em $\times_A$-bialgebra} if $(H, \Delta, \mu)$ is a $\times_A$-coalgebra
and $\Delta$ and $\mu$ are maps of $A^e$-rings.
\end{definition}

\section{$\times_A$-bialgebras versus bialgebroids}
The aim of this section is to clarify that the three notions
recalled in the previous section are in fact equivalent to each other.

\begin{theorem}\label{3.1}
For an $A^e$-ring $(H, s,t)$, the following data are equivalent :

(1) A bialgebroid structure $(H,s,t, \Delta, \eps)$;

(2) A bialgebroid with an anchor structure $(H,s,t, \Delta, \mu)$;

(3) A $\times_A$-bialgebra structure $(H, s, t, \Delta, \mu)$;

(4) A monoidal structure on ${}_H\M$ such that the forgetful functor
$F: {}_H\M \to {}_A\M_A$ is strict monoidal.
\end{theorem}

\begin{proof}
The equivalence $(3) \Leftrightarrow (4)$ is proven in
\cite[Theorem 5.1]{Sch:bia}.

$(1) \Rightarrow (2), (3)$. Let $(H,s,t, \Delta, \eps)$
be an $A$-bialgebroid in
the sense of Definition~\ref{def.Lu}, and define (cf.\
\cite[Eq.~(8)]{Lu})
\begin{equation}\label{37}
\mu= \mu_{\eps} : H\to \End (A), \quad
\mu (h) (a) = h \triangleright a
:= \eps (hs(a)) = \eps (ht(a)).
\end{equation}
The map $\mu$ is an algebra morphism since  $(A, \triangleright)
\in {}_H\M$.
The fact that $\mu$ is  $A$-bilinear follows by an elementary
calculation. Explicitly, for any
$a$, $b\in A$, $h\in H$ we have
$$
\mu (a\cdot h)(b) = \mu (s(a)h)(b) = \eps (s(a) h s(b))
 =  a \eps (h s(b)) = a\mu (h)(b)
=(a \cdot \mu(h))(b),
$$
where we used (\ref{8}) to derive the last equality, thus proving
that $\mu$ is left $A$-linear. Similar calculation that uses
(\ref{9'}), proves the right $A$-linearity of $\mu$.
Next we prove  that (A1) and (A2)
hold for $\mu$. Using (B2) and (\ref{4.2}) we have
$
\Delta (h s(a)) =\Delta (h) \Delta (s(a)) = h_{(1)} s(a) \otimes_A
h_{(2)}$.
Now using the first part of the counit property (\ref{24}) for $hs(a)$,
we obtain
$s\bigl( \eps ( h_{(1)} s(a) ) \bigl) h_{(2)} = h s(a)$,
i.e., (A1) for $\mu$. The condition (A2 ) follows from
 $h \triangleright a = \eps (ht(a))$ and  the second
part of the counit property  (\ref{24}) together with (B2) and (\ref{4.2}).
This shows that $(H,s,t, \Delta, \mu =\mu_{\eps})$ is an
$A$-bialgebroid with an anchor in the sense of Definition~\ref{2.4},
i.e.,  $(1) \Rightarrow (2)$.

In fact there is more, and this is $(1) \Rightarrow (3)$. Since
$\mu$ and
$\Delta$ are $A$-bimodule maps, they are left $A^e$-module maps.
Furthermore, both $\mu$ and the corestriction
$\Delta'$ of
$\Delta$ to
$\Gamma = H\times_A H$, are $k$-algebra maps. Therefore, by the
observation at the end of Section~\ref{subs.notation},
$\mu$ and $\Delta'$ are  maps of $A^e$-rings, and hence also maps of
$A^{e}$-bimodules.
Remark \ref{2.7} then implies that
$(H,s,t, \Delta', \mu =\mu_{\eps})$ is a $\times_A$-bialgebra.

$(2) \Rightarrow (1)$. Let  $(H,s,t, \Delta, \mu)$ be a
bialgebroid with an anchor and define
$\eps = \eps_{\mu} :H\to A$,
$h\mapsto \mu (h)(1_A)$.
\cite[Proposition 3.3]{Xu} shows that $\eps$ is a map in
${}_A\M_A$, and a
counit for $\Delta$.  This also implies that $\eps(s(a)) = \eps(t(a)) =a$,
hence, in particular that $\eps (1_H) =1_A$. Furthermore by \cite[Proposition
3.5]{Xu}, for all $g,h\in H$, $\eps(gh) = \mu(g)(\eps(h))$. Therefore $\eps(gh)
= \mu(g)(\eps(h)) = \mu(g)(\eps(s(\eps(h)))) = \eps(hs(\eps(h)))$, and similarly
for the target map $t$. This proves equations
(\ref{25}), and we conclude that
$(H,s,t, \Delta, \eps_{\mu})$ is a
bialgebroid in the sense of Definition \ref{def.Lu}.

The above result gives $\mu$ in terms of
$\eps = \eps_{\mu}$.  Looking at (\ref{37}) one can
 apply $(1) \Rightarrow (2)$ once again. In this way one obtains that
 for a
bialgebroid
with an anchor $(H,s,t, \Delta, \mu)$ the structure maps $\Delta'$ (the
corestriction
of $\Delta$ to $H\times_A H$) and $\mu$ are in fact maps of $A^e$-rings.
Following  Remark \ref{2.7} the implication $(3) \Rightarrow (2)$ is
then obvious.
\end{proof}

\section{Braided commutative algebras and bialgebroids}
Let $H$ be a bialgebra, $(A, \cdot)$ a left $H$-module algebra and let
$A\# H=A\otimes H$ as a $k$-module with the multiplication
$$
(a\# g) (b\# h)= a(g_{(1)}\cdot b) \# g_{(2)}h.
$$
Let ${}_H{\cal YD}^H$ denote the  pre-braided monoidal category
of (left-right) crossed or Yetter-Drinfeld modules. This means that
$(M, \cdot, \rho^M) \in {}_H{\cal YD}^H$ if and only if
$(M, \cdot)$ is a left
$H$-module, $(M, \rho^M)$ is a right $H$-comodule and
\begin{equation}
h_{(1)}\cdot m_{<0>} \otimes h_{(2)}m_{<1>} =
(h_{(2)} \cdot m)_{<0>} \otimes (h_{(2)} \cdot m)_{<1>} h_{(1)},
\label{yet}
\end{equation}
for all $h\in H$, $m\in M$.  For all $M$, $N\in {}_H{\cal YD}^H$,
$M\otimes N \in {}_H{\cal YD}^H$ via
$$
h\cdot (m\otimes n) = h_{(1)}\cdot m \otimes h_{(2)} \cdot n, \quad
m\otimes n \to m_{<0>} \otimes n_{<0>} \otimes n_{<1>} m_{<1>}.
$$
The pre-braiding (a braiding if $H$ has an antipode) is given by
$$
\sigma_{M, N}: M\otimes N\to N\otimes M, \quad
\sigma_{M, N}(m\otimes n) = n_{<0>} \otimes n_{<1>}\cdot m.
$$
An algebra $A$ which is also an object in ${}_H{\cal YD}^H$ via
$(A, \cdot, \rho^A)$, is called an
{\em algebra in ${}_H{\cal YD}^H$} if the algebra structures are maps in
the category ${}_H{\cal YD}^H$.
This is equivalent to say that $(A, \cdot)$ is
a left $H$-module algebra and $(A, \rho^A)$ is a right $H^{op}$-comodule
algebra.
An algebra $A$ in ${}_H{\cal YD}^H$ is said to be {\em braided commutative}
if the
multiplication $m_{A}$ of $A$ is commutative with respect to $\sigma_{A, A}$,
i.e.,
$m_{A}\circ \sigma_{A, A} = m_{A}$, or equivalently, for all $a,b\in A$,
\begin{equation}\label{45}
b_{<0>}(b_{<1>}\cdot a) = ab.
\end{equation}

The following theorem is a generalisation  of \cite[Theorem
5.1]{Lu}.

\begin{theorem}\label{4.1}
Let $H$ be a bialgebra, $(A, \cdot)$ a left $H$-module algebra and
$(A, \rho^A)$ a right $H$-comodule. Then  $(A, \cdot, \rho^A)$ is a
braided commutative algebra in ${}_H{\cal YD}^H$ if and only if
$(A\# H, s, t, \Delta, \eps)$ is an $A$-bialgebroid with the
source, target, comultiplication and the counit given by
$ s(a)= a\# 1_H$, $t(a)= a_{<0>}\# a_{<1>}$,
$\Delta (a\# h) = a\# h_{(1)} \otimes_A 1_A \# h_{(2)}$, $
\eps (a\# h)= \eps_H(h) a$,
for all $a\in A$ and $h\in H$.

Furthermore, if $H$ has an antipode $S$, then $A\# H$ is a Hopf
algebroid with the antipode
$$
\tau: A\# H \to A\# H, \quad
\tau (a\# h) = \Bigl (S(h_{(2)}) S^{2} (a_{<1>}) \Bigl) \cdot a_{<0>} \#
S(h_{(1)}) S^{2} (a_{<2>})
$$
for all $a \in A$ and $h\in H$.
\end{theorem}

\begin{proof}
Clearly, $s$ is an algebra map. We  prove now that $t$ is an anti-algebra map
if and only if $(A, \rho^A)$ is a right $H^{op}$-comodule algebra and the
braided
commutativity relation (\ref{45}) holds. Take any $a,b\in A$, then $
t(ab) = (ab)_{<0>} \# (ab)_{<1>}$, and
$t(b)t(a) = b_{<0>}(b_{<1>}\cdot a_{<0>}) \# b_{<2>}a_{<1>}$.
Suppose $t$ is an anti-algebra map. Then applying
$A\otimes \eps_H$ to the above equality one obtains  equation
(\ref{45}). It follows then that
$t(b) t(a) = a_{<0>}b_{<0>} \# b_{<1>}a_{<1>}$, i.e.
$\rho^A :A\to A\otimes H^{op}$ is an algebra map, hence $(A, \rho^A)$ is an
$H^{op}$-comodule algebra as required.
Conversely, suppose that $(A, \rho^A)$ is a right $H^{op}$-comodule algebra and
equation (\ref{45}) holds. Then
$$
t(b)t(a) = b_{<0>}(b_{<1>}\cdot a_{<0>}) \# b_{<2>}a_{<1>}=
a_{<0>}b_{<0>} \# b_{<1>}a_{<1>}= t(ab).
$$
Assume now that $t$ is an anti-algebra map.
We  prove that
${\rm Im}(\Delta) \subseteq \Gamma$ if and only if
$(A, \cdot, \rho^A) \in {}_H{\cal YD}^H$.
$A\# H$ is a right $A$-module via (\ref{16'}), i.e., using (\ref{45})
we have
\begin{equation}\label{46}
(b\# h)\cdot a = t(a) (b\# h) = a_{<0>} (a_{<1>}\cdot b) \# a_{<2>} h
= b a_{<0>} \# a_{<1>} h.
\end{equation}
${\rm Im}(\Delta) \subseteq \Gamma$ if and only if for all $a, b\in
A$,
$h\in H$ we have that
$(a\# h_{(1)}) (b_{<0>}\# b_{<1>} )\otimes_A 1\# h_{(2)} =
a\# h_{(1)} \otimes_A (1\# h_{(2)}) (b\# 1)$
or, equivalently,
$$
a(h_{(1)} \cdot b_{<0>})\# h_{(2)} b_{<1>} \otimes_A 1\# h_{(3)} =
a\# h_{(1)} \otimes_A h_{(2)}\cdot b\# h_{(3)}.
$$
Since the tensor product is defined over $A$
and equality (\ref{46}) holds we have
$$
a(h_{(1)} \cdot b_{<0>})\# h_{(2)} b_{<1>} \otimes_A 1\# h_{(3)} =
a (h_{(2)}\cdot b)_{<0>} \# (h_{(2)}\cdot b)_{<1>} h_{(1)} \otimes_A 1\#
h_{(3)}.
$$
Thus we conclude that  ${\rm Im}(\Delta) \subseteq \Gamma$ if and only if
$h_{(1)} \cdot b_{<0>}\otimes h_{(2)} b_{<1>} =
(h_{(2)}\cdot b)_{<0>} \otimes (h_{(2)}\cdot b)_{<1>} h_{(1)}$,
i.e., if and only if $(A, \cdot, \rho^A) \in {}_H{\cal YD}^H$.

Therefore we have proven that $t$ is an anti-algebra map and
${\rm Im}(\Delta) \subseteq \Gamma$ if and only if $(A, \cdot , \rho^A)$
is a braided commutative
algebra in ${}_H{\cal YD}^H$.
It is then straightforward to check that all the remaining conditions in
Definition~\ref{def.Lu}  hold.


Finally we  prove that $\tau$ defined in the theorem is the antipode of
$A\# H$. The canonical projection
$(A\# H) \otimes (A\# H) \to (A\# H) \otimes_A (A\# H)$ has a well defined
section
$\gamma : (A\# H) \otimes_A (A\# H) \to (A\# H) \otimes (A\# H)$,
$\gamma (a\# h \otimes_A b\# g ) = ab_{<0>} \# b_{<1>}h \otimes 1_A \#g$.
Since for all $a\in A$, $h\in H$,
$\tau (1_A\# h) = 1_A\# S(h)$ and
$\tau (a\# 1_H) = S^{2} (a_{<1>}) \cdot a_{<0>} \# S^{2} (a_{<2>})$,
we have
$$
\tau (a\# h) = \tau ((a\# 1_H)(1_A\# h) ) = \tau (1_A\# h) \tau (a\#
1_H),
$$
i.e., $\tau$ is an anti-algebra map.
Condition (ANT3) follows from the definition of $\gamma$,
$$
(a\# h_{(1)}) \tau (1_A \# h_{(2)}) =
(a\# h_{(1)}) (1_A \# S(h_{(2)})) = \eps_H (h) a\# 1_H,
$$
while (ANT1) can be established by  the following computation
\begin{eqnarray*}
\tau (t(a)) &=&\tau (a_{<0>} \# a_{<1>}) =
S(a_{<4>}) S^2(a_{<1>}) \cdot a_{<0>} \# S(a_{<3>}) S^2(a_{<2>}) \\
&=& S(a_{<4>}) S^2(a_{<1>}) \cdot a_{<0>} \# S\Bigl(S(a_{<2>}) a_{<3>}\bigl)\\
&=& S\Bigl(S(a_{<1>}) a_{<2>} \Bigl) \cdot a_{<0>} \# 1_H
= a\# 1_H = s(a).
\end{eqnarray*}
It remains to prove property (ANT2). The left hand side of (ANT2)
equals
$$
\tau (a\# h_{(1)}) (1_A \# h_{(2)}) =
S(h_{(2)}) S^2 (a_{<1>}) \cdot a_{<0>} \# S(h_{(1)}) S^2 (a_{<2>}) h_{(3)}.
$$
Equation~(\ref{yet}), evaluated at $\Delta(S(h\sw 1))\otimes h\sw 2 =
S(h\sw 2)\otimes S(h\sw 1)\otimes h\sw 3$ implies that
\begin{equation}\label{dinamo}
Sh_{(2)} \cdot a_{<0>} \otimes Sh_{(1)} a_{<1>}
h_{(3)} = \Bigl( S(h)\cdot a_{<0>}\Bigl)_{<0>}\otimes \Bigl( S(h)\cdot
a_{<0>}\Bigl)_{<1>}
\end{equation}
for all $h\in H$ and $a\in A$. Now
the right hand side of (ANT2) reads
\begin{eqnarray*}
(t \circ \eps \circ \tau) (a\# h) &=&
t(S(h) S^2(a_{<1>}) \cdot a_{<0>}) \\
&=& \Bigl(S(S(a_{<1>})h) \cdot a_{<0>} \Bigl)_{<0>} \#
\Bigl(S( S(a_{<1>})h) \cdot a_{<0>} \Bigl)_{<1>}\\
~~(\ref{dinamo})~~
&=& S((S(a_{<1>})h)\sw 2)\cdot a_{<0>} \#
  S((S(a_{<1>})h)\sw 1)a_{<1>}(S(a_{<1>})h)\sw 3\\
&=& S(h_{(2)}) S^2(a_{<3>}) \cdot a_{<0>} \#
S(h_{(1)}) S^2(a_{<4>}) a_{<1>}
S(a_{<2>}) h_{(3)}\\
&=& S(h_{(2)}) S^2(a_{<1>}) \cdot a_{<0>} \#
S(h_{(1)}) S^2(a_{<2>}) h_{(3)}
\end{eqnarray*}
that is exactly the left hand side of (ANT2). Hence, $\tau$ is an
antipode of $A\# H$.
\end{proof}

Theorem~\ref{4.1} generalises,
gives  converse to, and a more transparent proof of \cite[Theorem 5.1]{Lu}.
It also
provides one with a rich source of examples of bialgebroids.
Several  examples of braided commutative algebras in
${}_H{\cal YD}^H$ are known, cf.\
\cite{COZ}, \cite{CW}, \cite{Maj90}.
For example, for an $H^{op}$-Galois
extension,
$A/B$, the centralizer algebra $E=C_A(B)$ has a structure of a braided
commutative
algebra in ${}_H{\cal YD}^H$ (\cite{COZ}).

We indicate now three other ways of obtaining braided commutative
algebras.

\begin{example}\label{steaua}
1. Let $(H, R= \sum R^1\otimes R^{2})$ be a quasitriangular
bialgebra and $(A, \cdot)$ a left
$H$-module algebra, which is braided commutative in the
pre-braided category ${}_H\M$: i.e.
$\sum (R^2 \cdot b)(R^1\cdot a)=ab$, for all $a$, $b\in A$.
Then $A$ is a braided commutative algebra in
${}_H{\cal YD}^H$ where the coaction of $H$ on $A$ is given by
$\rho^{A} :A \to A\otimes H$, $a\mapsto \sum R^2\cdot a \otimes R^1$.
In this way, all examples of braided commutative
algebras over a quasi-triangular Hopf algebra $H$
from \cite{CW} give
examples of braided commutative algebras in  ${}_H{\cal YD}^H$, and
hence examples of bialgebroids $A\#H$.
Note that \cite[Theorem 5.1]{Lu}
corresponds to the quasi-triangular Hopf algebra
$(D(H), {\cal R})$.

2. Dually, let $(H, \sigma)$ be a coquasitriangular bialgebra and
$(A, \rho^{A})$ be a right $H^{op}$-comodule algebra such that
$ab = \sigma (a_{<1>} \otimes b_{<1>}) b_{<0>}a_{<0>}$,
for all $a$, $b\in A$.
Then $(A, \cdot, \rho^{A})$ is a braided commutative algebra in
${}_H{\cal YD}^H$, where the left $H$-action is
$h\cdot a = \sigma (a_{<1>} \otimes h) a_{<0>}$.

3. There is a general way of constructing braided commutative algebras
in ${}_H{\cal YD}^H$ pointed out in \cite{Maj90}, \cite{COZ}.
Let $(V, \cdot, \rho^V) \in {}_H{\cal YD}^H$ and $T(V)$ be
the tensor algebra of $V$. Then the (co)-actions of $H$ on $V$
extend uniquely to (co)-actions on $T(V)$ such that $T(V)$ becomes
an algebra in the category ${}_H{\cal YD}^H$. Let $S^b(V)$
be the ``braided symmetric" algebra of $V$, i.e., $S^b(V): = T(V)/I$,
where $I$ is the two-sided ideal of $T(V)$ generated by all
elements of the form $
v\otimes w - w_{<0>}\otimes w_{<1>} \cdot v$,
for all $v, w \in V$. Then $S^b(V)$ is a braided commutative
algebra in ${}_H{\cal YD}^H$.
\end{example}

Using the  FRT-construction and Example \ref{steaua} we
present now a generic construction of bialgebroids associated
to any solution of the quantum Yang-Baxter equation.

Let $n$ be a positive integer and
$R=(R_{uv}^{ij}) \in M_n(k) \otimes M_n(k)$ be a solution of the
QYBE, $R^{12}R^{13}R^{23}=R^{23}R^{13}R^{12}$. Let $A(R)$ be the
bialgebra associated to $R$ using the FRT construction:
$A(R)$ is a free $k$-algebra generated be
$(c^i_j)_{i,j=1, \cdots, n}$ with the relations
$R_{vu}^{ij} c_k^u c_l^v = R_{lk}^{vu} c_v^i c_u^j$,
for all $i$, $j$, $k$, $l=1, \cdots, n$ (Einstein's summation
convention assumed), and the standard matrix bialgebra structure.
View $R$ as an endomorphism
of $V\otimes V$, for an $n$-dimensional vector space, and define the
corresponding braided symmetric algebra $S^b(V) =S_{R}(n)$  as follows.
$S_R(n)$ is a free $k$-algebra generated by
$\xi_1, \cdots, \xi_n$ with the relations
$\xi_u \xi_v =R_{uv}^{li}\xi_i\xi_l$,
for all $u$, $v=1, \cdots, n$. $S_R(n)$ is a braided commutative algebra in
${}_{A(R)}{\cal YD}^{A(R)}$ via
$c_v^j \cdot \xi_u = R_{uv}^{ij} \xi_i$ and $\xi_v \to \xi_u \otimes c_v^u
$, for all  $j$, $u$, $v=1, \cdots, n$. Then  Theorem~\ref{4.1} implies
the following example of a quantum groupoid

\begin{proposition}\label{parma}
Let $n$ be a positive integer and $R=(R_{uv}^{ij})$ a solution
of the QYBE. Then the smash product $S_R(n) \# A(R)$ has
a structure of $S_R(n)$-bialgebroid with the
source, target, comultiplication and the counit given by
$s(\xi_i)= \xi_i\# 1$, $t(\xi_i)= \xi_u \# c^u_i$,
$\Delta (\xi_i\# c^u_v) = \xi_i\# c^u_l \otimes_{S_R(n)} 1 \# c^l_v$, and
$\eps (\xi_i \# c^u_v)= \delta_{u,v} \xi_i$, where
$i, u, v=1, \cdots, n$.
\end{proposition}

In particular, Proposition \ref{parma} associates bialgebroids to quantum
matrix groups such as $GL_{q}(n)$ and their corepresentation spaces
such as the quantum hyperplane (cf. \cite{Maj:book}).

\section{Comments on duals of bialgebroids - bicoalgebroids}
On formal level, the notion of a bialgebra is self-dual in the
following sense. Write definition of a bialgebra in terms of
commutative diagrams.
Then the structure obtained by reversing arrows in diagrams defining a
bialgebra is again a bialgebra.  It is clear that, in general,
a dual of a bialgebroid in the above sense is no longer a
bialgebroid.\footnote{This should not be confused with a left or right 
dual module of a bialgebroid which is a bialgebroid provided certain 
finitely generated projective type conditions are satisfied (cf.\ 
\cite{KadSzl:dua}).} This is because by reversing the arrows in diagrams 
defining an
algebra and a module one obtains diagrams defining a coalgebra and a 
comodule. Thus
if one wants to construct a (formally) dual object to a bialgebroid one has to
consider an object within the category of comodules of a coalgebra.

\begin{definition}
    Let $C$ be a coalgebra over a field $k$. A {\em bicoalgebroid} is
    a $k$-coalgebra $H$ which satisfies the following conditions:
    \begin{enumerate}
        \item[(BC1)] There is a coalgebra map $\alpha :H\to C$ and an
        anti-coalgebra
        map $\beta:H\to C$ such that for all $h\in H$,
        $
        \alpha(h\sw 1)\otimes \beta(h\sw 2) = \alpha(h\sw 2)\otimes
        \beta(h\sw 1).
        $
        This allows one to view $H$ as a $C$-bicomodule via
        left  coaction
        $
        {}^{H}\rho(h) = \alpha(h\sw 1)\otimes h\sw 2$, and the right coaction
        $\rho^{H}(h) =
        h\sw 2\otimes \beta(h\sw 1)$.
        Let
        $$
        H\square_{C}H =\{ \sum_{i}g^{i}\otimes h^{i}\in H\otimes H \; |\;
\sum_{i}g^{i}\sw
        2\otimes \beta(g^{i}\sw 1)\otimes h^{i} = \sum_{i}g^{i}\otimes
        \alpha(h^{i}\sw 1)\otimes h^{i}\sw 2\}
        $$
        be the corresponding cotensor product.
        \item[(BC2)] There is a $C$-bicomodule map $\mu :H\square_{C}H\to H$ which
        is an associative product with respect to the cotensor product and
such
        that for all $\sum_{i}g^{i}\otimes h^{i}\in H\square_{C}H$:
        \begin{enumerate}
            \item $
        \sum_{i}\mu(g^{i}\otimes h^{i}\sw 1)\otimes \alpha(h^{i}\sw 2) =
        \sum_{i}\mu(g^{i}\sw 1\otimes h^{i})\otimes \beta(g^{i}\sw 2),
        $
            \item
        $
        \Delta(\mu(\sum_{i}g^{i}\otimes h^{i})) = \sum_{i}\mu(g^{i}\sw
        1\otimes h^{i}\sw 1)\otimes \mu(g^{i}\sw 2\otimes h^{i}\sw 2).
        $
        \end{enumerate}
        \item[(BC3)] There exists a bicomodule map $\eta :C\to H$ which is a unit
        for $\mu$, i.e.,
        $$
        \mu\circ(\eta\square_{C}H)\circ{}^{H}\rho =
        \mu\circ(H\square_{C}\eta)\circ\rho^{H} = H,
        $$
        and such that for all $c\in C$, $\eps(\eta(c)) =\eps(c)$, and
        $$
        \Delta(\eta(c)) = \eta(c)\sw 1\otimes \eta(\alpha(\eta(c)\sw 2)) =
        \eta(c)\sw 1\otimes \eta(\beta(\eta(c)\sw 2)).
        $$
    \end{enumerate}
\label{def.bicoalgebroid}
\end{definition}

Few comments are needed in order to see that the above definition
makes sense. The condition (BC2)(a) makes sense because for all
$\sum_{i}g^{i}\otimes h^{i}\in H\square_{C}H$ we have
\begin{equation}
\sum_{i} g^{i}\sw 2\otimes g^{i}\sw 3 \otimes \beta(g^{i}\sw 1)\otimes
h^{i} = \sum_{i}g^{i}\sw 1\otimes g^{i}\sw 2\otimes \alpha(h^{i}\sw
1)\otimes h^{i}\sw 2,
\label{D1}
\end{equation}
\begin{equation}
\sum_{i} g^{i}\sw 2\otimes\beta(g^{i}\sw 1)\otimes h^{i}\sw 1\otimes
h^{i}\sw 2 = \sum_{i} g^{i}\otimes \alpha(h^{i}\sw 1)\otimes h^{i}\sw
2\otimes h^{i}\sw 3
\label{D2}
\end{equation}
Equation (\ref{D1}) implies that $\sum_{i} g^{i}\sw 1\otimes
h^{i}\otimes g^{i}\sw 2 \in H\square_{C}H\otimes H$ while Equation (\ref{D2})
implies that $\sum_{i}g^{i}\otimes h^{i}\sw 1\otimes h^{i}\sw 2 \in
H\square_{C}H\otimes H$. Furthermore both equations (\ref{D1}) and
(\ref{D2}) imply that $\sum_{i}g^{i}\sw 1\otimes h^{i}\sw 1\otimes
g^{i}\sw 2\otimes h^{i}\sw 2 \in H\square_{C}H\otimes H\otimes H$.
Using condition (BC2)(a) one concludes that $\sum_{i}\mu(g^{i}\sw
1\otimes h^{i}\sw 1)\otimes g^{i}\sw 2\otimes h^{i}\sw 2 \in H\otimes
H\square_{C}H$, i.e., condition (BC2)(b) makes sense.
Note that conditions (BC2) and (BC3)
mean  also that $H$ is a $C$-ring
in the sense of \cite[Section~6]{Brz:cor}.

One way of understanding the relation of an object defined in
Definition~\ref{def.bicoalgebroid} to bialgebroids is to
write all the conditions in terms of commutative diagrams. Reversing
the arrows, replacing $C$ by $A$, $\alpha$ by $s$, $\beta$ by $t$,
$\square_{C}$ by $\otimes_{A}$,
$\mu$ by $\Delta$ and $\eta$ by $\epsilon$ one obtains commutative
diagrams defining a bialgebroid.

An indication that this dualisation of a bialgebroid might play 
a role in introducing  self-dual bialgebroids
comes from the
following observation. A self-dual generalisation of a Hopf algebra is
provided by the notion of a {\em weak Hopf algebra}
\cite{Boh:wea}. A {\em weak bialgebra} is a vector space $H$ which is an
algebra and a coalgebra with multiplicative (but non-unital) coproduct such
that for
all $x,y,z\in H$,
$\eps(xyz) = \eps(x y\sw 1)\eps(y\sw 2 z) =
\eps(x y\sw 2)\eps(y\sw 1 z)$,
and
$(\Delta\otimes H)\circ \Delta(1) = (\Delta(1)\otimes 1)(1\otimes
\Delta(1)) = (1\otimes \Delta(1))(\Delta(1)\otimes 1)$.
A weak bialgebra $H$ is a {\em weak Hopf algebra} if there exists an
antipode, i.e., a
linear map $S:A\to A$ such that for all $h\in H$,
$
h\sw 1 S(h\sw 2) = \eps(1\sw 1h)1\sw 2$, $S(h\sw 1)h\sw 2 = 1\sw
1\eps(h 1\sw 2)$, and $S(h\sw 1) h\sw 2 S(h\sw 3) = S(h)$.
Weak Hopf algebras have been studied in connection to
integrable models and classification of subfactors of von Neumann
algebras. In \cite[Proposition~2.3.1]{EtiNik:dyn}
it has been shown that a weak Hopf algebra
with bijective antipode is a bialgebroid over $A = {\rm Im}\eps_{t}$,
where $\eps_{t}: H\to H$, $h\mapsto \eps(1\sw 1h)1\sw 2$.
By \cite[Eq.\ (2.12)]{Boh:wea} $\ker\eps_{t}$ is a coideal,
hence we can state
the following

\begin{proposition}
    Let $H$ be a weak Hopf algebra with bijective antipode. Let $C =
    H/\ker\eps_{t}$ with the canonical surjection $\pi_{t}:H \to C$.
    Then $C$ is a coalgebra and
    $H$ is a bicoalgebroid over $C$ with the following structure maps:
    \begin{enumerate}
        \item $\alpha = \pi_{t}$, $\beta = \pi_{t}\circ S^{-1}$,
        \item $\mu (\sum_{i}g^{i}\otimes h^{i}) = \sum_{i}g^{i}h^{i}$, for
        all $\sum_{i}g^{i}\otimes h^{i}\in H\square_{C}H$.
        \item $\eta :C\to H$, $c\mapsto \eps_{t}(h)$, where $h\in
        \pi^{-1}_{t}(c)$.
     \end{enumerate}
\end{proposition}
\begin{proof}
    This can be proven by dualising the proof of
\cite[Proposition~2.3.1]{EtiNik:dyn}. Although not elementary this is quite
straightforward and we leave it to the reader\footnote{The authors
will be
happy to supply the full proof upon request.}.
\end{proof}

Thus  a weak Hopf algebra is an example of both a bialgebroid and a
bicoalgebroid.
This suggests that if one imposes a selfduality as a key property that
must be enjoyed by a proper generalisation of
a bialgebra, such a generalisation
should be a bialgebroid over an algebra $A$ and  a
bicoalgebroid over a coalgebra $C$ at the same time. Some relations
between $A$ and $C$ should also be required in order to compare tensor
products with cotensor products.
Once such
a relationship is imposed compatibility conditions between product and
coproduct must involve both  tensor and cotensor products. What
these should be we consider an interesting open question.

\begin{center}
{\sc Acknowledgements}
\end{center}
Tomasz Brzezi\'nski would like to thank
EPSRC for an Advanced Research Fellowship, and Gigel Militaru would
like to thank the Royal Society for a visiting fellowship.

\end{document}